\documentclass[12pt,reqno]{article}

\usepackage[usenames]{color}
\usepackage{amssymb}
\usepackage{graphicx}
\usepackage{amscd}

\usepackage{amsthm}
\newtheorem{theorem}{Theorem}

\newtheorem{proposition}[theorem]{Proposition}

\theoremstyle{definition}

\newtheorem{example}[theorem]{Example}

\usepackage[colorlinks=true,
linkcolor=webgreen, filecolor=webbrown,
citecolor=webgreen]{hyperref}

\definecolor{webgreen}{rgb}{0,.5,0}
\definecolor{webbrown}{rgb}{.6,0,0}

\usepackage{color}

\usepackage{float}

\usepackage{graphics,amsmath,amssymb}
\usepackage{amsfonts}
\usepackage{latexsym}
\usepackage{epsf}

\newcommand{\seqnum}[1]{\href{http://oeis.org/#1}{\underline{#1}}}

\begin{document}

\begin{center}
\vskip 1cm{\LARGE\bf  On the Hankel transform of C-fractions}  \vskip 1cm
\large
Paul Barry\\
School of Science\\
Waterford Institute of Technology\\
Ireland\\

\end{center}
\vskip .2 in

\begin{abstract} We study the Hankel transforms of sequences whose generating function can be expressed as a C-fraction. In particular, we relate the index sequence of the non-zero terms of the Hankel transform to the powers appearing in the monomials defining the C-fraction. A closed formula for the Hankel transforms studied is given. As every power-series can be represented by a C-fraction, this gives in theory a closed form formula for the Hankel transform of any sequence. The notion of multiplicity is introduced to differentiate between Hankel transforms.
\end{abstract}

\section{Introduction} Given a sequence $a_n$, we denote by $h_n$ the general term of the sequence
with $h_n=|a_{i+j}|_{0 \le i,j \le n}$. The sequence $h_n$ is called the Hankel transform of $a_n$ \cite{Kratt_1, Kratt_2, Layman}. If the sequence $a_n$ has generating function $g(x)$, then by an abuse of language we can also refer to $h_n$ as the Hankel transform of $g(x)$.

A well known example of Hankel transform is that of the Catalan numbers, $C_n=\frac{1}{n+1}\binom{2n}{n}$, where we find that
$h_n=1$ for all $n$. Hankel determinants occur naturally in many branches of mathematics, from combinatorics \cite{Brualdi} to number theory \cite{Milne} and to mathematical physics \cite{Vein_Dale}.

We shall be interested in characterizing the Hankel transform of sequences whose generating functions can be expressed as the following type of C-fraction:
\begin{equation}\label{Cfract1} g(x)=\cfrac{1}
{1+\cfrac{a_1 x^{q_1}}{1+\cfrac{a_2 x^{q_2}}{1+\cfrac{a_3 x^{q_3}}{1+\cdots}}}},\end{equation} for appropriate values of coefficients $a_1,a_2,a_3,\ldots$ and exponents $q_1,q_2,q_3,\ldots$.
The results will depend on making explicit the relationship between this type of C-fraction, and $h(1/x)$, where $h(x)$ is the following type of continued fraction:
\begin{equation}\label{Cfract2}
h(x)=\cfrac{x^{p_0}}
{b_1 x^{p_1}+\cfrac{1}{b_2 x^{p_2}+\cfrac{1}{b_3 x^{p_3}+\cfrac{1}{b_4 x^{p_4}+\cdots}}}}.\end{equation}  We will then be able to use classical results \cite{Holtz} to conclude our study and to examine interesting examples.
\section{Review of known results}
The first part of this section reviews the close link between power series and C-fractions. Note that the ``C'' comes from the word ``corresponding''.

We commence with a power series
\begin{equation} \label{PS} f_0(x)=1+c_1 x+c_2 x^2+c_3 x^3+\ldots. \end{equation}
We form the family of power series $\{f_n(x)\}$ by the relations
\begin{equation} \label{Pn}
f_{n+1}(x)=\frac{a_{n+1}x^{q_{n+1}}}{f_n(x)-1}, \quad n=0,1,2,\ldots,\end{equation} where the $q_n$ are positive integers chosen together with complex numbers $a_n$ in such a way that if $f_n(x)\ne 1$, $f_{n+1}(0)=1$. If no $f_n(x)=1$, this process yields an infinite sequence of power series $f_0(x), f_1(x), f_2(x),\ldots$. If some $f_n(x)=1$, the process terminates and yields a finite set of power series $f_0(x), f_1(x), \ldots, f_n(x)$. The continued fraction
\begin{equation} \label{Cfract3} 1+\cfrac{a_1 x^{q_1}}
{1+\cfrac{a_2 x^{q_2}}{1+\cfrac{a_3 x^{q_3}}{1+\cfrac{a_4 x^{q_4}}{1+\cdots}}}},\end{equation} formed with these $a_n$ and $q_n$ is said to \emph{correspond} to the power series (\ref{PS}) \cite{Jones, Leighton}.
Conversely, if we begin with a continued fraction of the form (\ref{Cfract3}), we can form the $n$-th \emph{approximant} $\frac{A_n(x)}{B_n(x)}$ by means of the recurrence relations
\begin{equation*} \begin{split}
A_0 =1, \quad \quad B_0=1,\\
A_1=1+a_1 x^{q_1},\quad \quad B_1=1,\\
A_n=A_{n-1}+a_nx^{q_n}A_{n-2},\quad B_n=B_{n-1}+a_n x^{q_n}B_{n-2},\\
\quad \quad\quad n=2,3,\ldots.\end{split} \end{equation*}
We have
\begin{equation}\label{DET}
\frac{A_n(x)}{B_n(x)}-\frac{A_{n-1}(x)}{B_{n-1}(x)}=\frac{(-1)^{n-1}a_1 a_2 a_3\cdots a_n x^{s_n}}{B_{n-1}(x)B_{n-2}(x)}, \end{equation}
where
$$s_n=q_1+q_2+\cdots+ q_n.$$ By equation (\ref{DET}) the Taylor development of the rational function $\frac{A_{n-1}(x)}{B_{n-1}(x)}$ about the origin agrees with the development of $\frac{A_n(x)}{B_n(x)}$ up to but not including the term in $x^{s_n}$. Hence if (\ref{Cfract3}) is nonterminating, the C-fraction (\ref{Cfract3}) determines uniquely a \emph{corresponding} power series.

We have the following classical result \cite{Leighton}
\begin{proposition} If the continued fraction (\ref{Cfract3}) corresponds to the power series (\ref{PS}), then the power series (\ref{PS}) corresponds to the continued fraction (\ref{Cfract3}), and conversely.
\end{proposition}
A division-free algorithm for the construction of the C-fraction (\ref{Cfract3}) from the power series (\ref{PS}) is given by Frank \cite{Frank_1, Frank_2}.

If we start with a power series $f(x)=\sum_{i=0}^{\infty}t_i x^i$, then by considering the sequence $1+xf(x)$, which is in the form (\ref{PS}), we see that $f(x)$ corresponds to a C-fraction of the form
$$\cfrac{a_0 x^{q_0}}{1+\cfrac{a_1 x^{q_1}}{1+\cfrac{a_2 x^{q_2}}{1+\cdots}}},$$ for appropriate values of
$a_0, a_1,a_2,a_3,\ldots$ and $q_0, q_1, q_2, q_3,\ldots$.

We now recall known results concerning the Hankel transform of sequences whose generating functions are of the form $f(1/x)$ where $f(x)$ can be expressed as a continued fraction of the form
\begin{equation}\label{Cfract4}
f(x)=\cfrac{b_0 x^{p_0}}
{b_1 x^{p_1}+\cfrac{1}{b_2 x^{p_2}+\cfrac{1}{b_3 x^{p_3}+\cfrac{1}{b_4 x^{p_4}+\cdots}}}}.\end{equation}
We have the following result \cite{Holtz}.
\begin{proposition}\label{Hankel} Let $h_n$ denote the Hankel transform of the sequence $[x^n]f(1/x)$ where $f(x)$ has the form (\ref{Cfract4}) (give conditions on $b_0=1$ and $p_0=0$).  Then $h_n$ is zero for all $n$ unless $n=p_1+p_2+\cdots+p_m$, for some $m$, in which case
\begin{equation} \label{H} h_n=\prod_{i=1}^m (-1)^{\frac{p_i(p_i-1)}{2}}\cdot (-1)^{\sum_{i=0}^{m-1} i p_{i+1}}\prod_{i=1}^m \frac{1}{b_i^{p_i+2 \sum_{j=i+1}^m p_j}}. \end{equation}
\end{proposition}
\section{Main result}
In order to obtain our main result, we need to relate C-fractions of the form
$$g(x)=\cfrac{1}
{1+\cfrac{a_1 x^{q_1}}{1+\cfrac{a_2 x^{q_2}}{1+\cfrac{a_3 x^{q_3}}{1+\cdots}}}}$$ to continued fractions of the form
$$f(x)=\cfrac{x^{p_0}}
{b_1 x^{p_1}+\cfrac{1}{b_2 x^{p_2}+\cfrac{1}{b_3 x^{p_3}+\cfrac{1}{b_4 x^{p_4}+\cdots}}}}.$$
We wish to find the conditions under which $f(1/x)=g(x)$. We look at the case of unit coefficients first. By equation (\ref{H}), the corresponding Hankel transforms will then take on values from the set $\{-1,0,1\}$.

By successive divisions above and below the line, we can cast $f(x)$ in the form
$$f(x)=\cfrac{x^{p_0-p_1}}{1+\cfrac{x^{-p_1-p_2}}{1+\cfrac{x^{-p_2-p_3}}{1+\cdots}}},$$  and hence we have
$$f(1/x)=\cfrac{x^{-p_0+p_1}}{1+\cfrac{x^{p_1+p_2}}{1+\cfrac{x^{p_2+p_3}}{1+\cdots}}}.$$ Starting from $g(x)$ and proceeding to $f(x)$  is more problematic, since it is not clear what to choose as $p_0$. The Hankel transforms that we will be concerned with determine that we require the condition $-p_0+p_1=0$, and hence that $p_1=p_0$. We choose to set $p_0=1$. Then starting from the C-fraction
$$\cfrac{1}
{1+\cfrac{ x^{q_1}}{1+\cfrac{ x^{q_2}}{1+\cfrac{ x^{q_3}}{1+\cdots}}}}$$ we find the following continued fraction of type (\ref{Cfract2}):
$$\cfrac{x^{p_0}}
{x^{p_0}+\cfrac{1}{x^{q_1-p_0}+\cfrac{1}{x^{q_2-q_1+p_0}+\cfrac{1}{x^{q_3-q_2+q_1-p_0}+\cfrac{1}{x^{q_4-q_3+q_2-q_1+p_0}+\cdots}}}}}$$
By Proposition (2), the position of the non-zero terms of the corresponding Hankel transform will be given by the indexing sequence $p_0, p_0+(q_1-p_0), p_0+(q_1-p_0)+(q_2-q_1+p_0), p_0+(q_1-p_0)+(q_2-q_1+p_0)+(q_3-q_2+q_1-p_0), \ldots$ or
$p_0, q_1, q_2+p_0, q_3+q_1, q_4+q_2+p_0, \ldots$. This sequence can be realised by
\begin{displaymath}\left(\begin{array}{c} m_0 \\ m_1 \\m_2\\ m_3 \\m_4 \\m_5 \\ \vdots \end{array}\right)=
\left(\begin{array}{ccccccc}1 & 0 & 0 & 0 & 0 & 0 & \ldots \\
                            0 & 1 & 0 & 0 & 0 & 0 & \ldots \\
                            1 & 0 & 1 & 0& 0 & 0 & \ldots \\
                            0 & 1 & 0 & 1 & 0 & 0 & \ldots \\
                            1 & 0 & 1&  0 & 1 & 0 & \ldots \\
                            0 & 1  & 0 & 1 & 0 & 1 &\ldots\\ \vdots
& \vdots &
\vdots & \vdots & \vdots & \vdots &
\ddots\end{array}\right)\left(\begin{array}{c} p_0 \\ q_1 \\q_2\\ q_3 \\q_4 \\q_5 \\ \vdots \end{array}\right).\end{displaymath}
The $n$-th term of this sequence $m_n$ is given by $$m_n=\sum_{k=0}^n \frac{1+(-1)^{n-k}}{2} \tilde{q}_k=\sum_{k=0}^n \sum_{i=0}^k (-1)^{k-i}\tilde{q}_i=\sum_{k=0}^n p_k,$$ where $\tilde{q}_0=p_0$,
$\tilde{q}_n=q_n$ for $n>0$, and $p_n=\sum_{k=0}^n (-1)^{n-k} \tilde{q}_k$.  Note that since the above matrix is $\left(\frac{1}{1-x^2},x\right)$ as a Riordan array, then if the g.f. of the sequence $q_1, q_2, q_3,\ldots$ is $G(x)$, then the g.f. of the index set is
$$\frac{1}{1-x^2}(1+xG(x)).$$
We next note that if
$$f(x)=\cfrac{b_0 x^{p_0}}
{b_1 x^{p_1}+\cfrac{1}{b_2 x^{p_2}+\cfrac{1}{b_3 x^{p_3}+\cfrac{1}{b_4 x^{p_4}+\cdots}}}}$$ is to be such that
$f(1/x)$ can be represented as
$$g(x)=\cfrac{a_0x^{q_0}}
{1+\cfrac{a_1 x^{q_1}}{1+\cfrac{a_2 x^{q_2}}{1+\cfrac{a_3 x^{q_3}}{1+\cdots}}}}$$  then we must have
\begin{equation}\label{ab} a_k = \frac{1}{b_k b_{k+1}}.\end{equation}
Reversing this set of equations, beginning with $b_0=1$, we find that
$$b_{2n}=\frac{a_0 a_2 \cdots a_{2n-2}}{a_1 a_3 \cdots a_{2n-1}},$$ and
$$b_{2n+1}=\frac{a_1 a_3 \cdots a_{2n-1}}{a_0 a_2 \cdots a_{2n}}.$$
(See also \cite{Krushchev}, Theorem $3.6$ and its corollaries).
 Substituting these values into Equation (\ref{H}) and simplifying (where we take $a_0=1$, $p_0=1$), gives us the main result of this note.
\begin{proposition} The non-zero elements of the Hankel transform of the sequence with generating function given by the C-fraction $$\cfrac{1}
{1+\cfrac{ a_1x^{q_1}}{1+\cfrac{a_2 x^{q_2}}{1+\cfrac{a_3 x^{q_3}}{1+\cdots}}}}$$ are given by
$$h_n=\prod_{i=1}^m (-1)^{\frac{p_i (p_i+1)}{2}}\cdot (-1)^{1+\sum_{i=0}^{m-1} i p_{i+1}}\cdot \prod_{k=1}^m a_k^{\sum_{i=k}^m p_i},$$ where
$$p_i=\sum_{j=0}^i (-1)^{i-j}\tilde{q}_j \quad \text {and} \quad n=\sum_{k=0}^m p_i,$$ and the sequence $\tilde{q_n}$ is given by ${1,q_1,q_2,q_3,\ldots}$.
\end{proposition}
\begin{example} We consider the Fibonacci-inspired C-fraction
$$\cfrac{1}
{1+\cfrac{ x}{1+\cfrac{ x}{1+\cfrac{ 2x^2}{1+\cfrac{3x^3}{1+\cdots}}}}}$$ where $\tilde{q}_n=F_n+0^n$ and
$a_n=F_n$. Then we find that the non-zero terms of the Hankel transform are indexed by
$$\sum_{k=0}^n {p_k}=\sum_{k=0}^n \sum_{i=0}^k (-1)^{k-i} (F_i+0^i)=F_{i+1}.$$ The non-zero terms, calculated as
$$\prod_{i=1}^m (-1)^{\frac{F_i (F_i+1)}{2}}\cdot (-1)^{1+\sum_{i=0}^{m-1} i F_{i+1}}\cdot \prod_{k=1}^m F_k^{\sum_{i=k}^m p_i},$$ begin
$$1, 1, 1, -2, 72, 1944000, \ldots.$$ We see emerging here an interesting phenomenon, which we will naturally term ``multiplicity''. In this case we note that $F_2=F_3=1$, corresponding to the first two $1$'s of the non-zero Hankel elements above. The Hankel transform is thus given by
$$1, 1, -2, 0, 72, 0, 0, 1944000, 0, 0, 0, 0, 1547934105600000000, 0, 0, \ldots,$$ where the initial $1$ has multiplicity $2$.
\end{example}
\begin{example} The Catalan numbers $C_n=\frac{1}{n+1}\binom{2n}{n}$ provide an interesting example of the notion of multiplicity. They have generating function
$$\cfrac{1}{1-\cfrac{x}{1-\cfrac{x}{1-\cdots}}},$$ and thus
$a_n=-1$ for all $n > 0$ and $q_n=1$ for all $n>0$ (and hence $\tilde{q}_n=1$ for all $n \ge 0$). We then have
$$\sum_{k=0}^n {p_k}=\sum_{k=0}^n \sum_{i=0}^n (-1)^{k-i} = {\lfloor \frac{n+2}{2} \rfloor}.$$ Thus the non-zero terms of the Hankel transform of $C_n$ are indexed by
$$1,1,2,2,3,3,4,4,5,5,\ldots$$
These terms are all calculated to equal $1$, as is well-known. Thus we can write the Hankel transform of $C_n$ as
$$1_2,1_2,1_2, \ldots,$$ where the sub-index $2$ indicates that each $1$ occurs with ``multiplicity'' $2$. This is a shorthand way of saying that the index set is $1,1,2,2,3,3,\ldots$.
\end{example}
\begin{example} It is well known that the Hankel transform of the aerated Catalan numbers $C_{\frac{n}{2}}\frac{1+(-1)^n}{2}$ is also the all-$1$'s sequence. This sequence has generating function
$$\cfrac{1}{1-\cfrac{x^2}{1-\cfrac{x^2}{1-\cdots}}},$$ where now
$a_n=-1$ for all $n > 0$ and $q_n=2$ for all $n>0$ (and hence $\tilde{q}_n=2-0^n$ for all $n \ge 0$). Now
$$\sum_{k=0}^n {p_k}=\sum_{k=0}^n \sum_{i=0}^n (-1)^{k-i} (2-0^i) = n+1, $$ and hence the indexing set for this Hankel transform is $1,2,3,4,5,6,\ldots$. That is, each $1$ appears with multiplicity $1$. Thus in a sense this is the original sequence with Hankel transform of all $1$'s.
\end{example}
\begin{example} The generalized Rogers-Ramanujan continued fraction. We consider the continued fraction
$$\cfrac{1}
{1+ \cfrac{\gamma x}{1+\cfrac{ \gamma x^2}{1+\cfrac{\gamma  x^3}{1+ \cfrac{\gamma x^4}{1+\cdots}}}}}.$$ Here,
$\tilde{q}_n=n+0^n$, and $a_n=\gamma-\gamma 0^n= \gamma (1-0^n)$.
We then have that $p_n$ is the sequence
$$1, 0, 2, 1, 3, 2, 4, 3, 5, 4, 6,\ldots,$$ and $\sum_{k=0}^n p_k $ is the sequence that begins
$$1, 1, 3, 4, 7, 9, 13, 16, 21, 25, 31, \ldots.$$ The non-zero terms of the Hankel transform are, in order,
$$1, 1, - \gamma^6, \gamma^{12}, \gamma^{32}, \gamma^{52}, - \gamma^{94}, \gamma^{136}, \gamma^{208}, \gamma^{280}, - \gamma^{390},\ldots.$$ The exponent sequence
$$0,0,6,12,32,52,94,\ldots$$ can be shown to have generating function
$$  \frac{2x^2(x^3+3)}{(x+1)^2(x-1)^4}.$$
\end{example}
\section{Conclusion}
Since to each sequence $a_n$ there corresponds the power series $\sum_{k=0}^{\infty} a_n x^n$, and to each power series there corresponds a C-fraction, the foregoing gives, in theory, a closed form formula for the Hankel transform of each sequence. Of course, this presupposes that the passage from generating function to C-fraction can be effected easily. The Q-D algorithm is one method for this.

We note that Heilermann's formula \cite{Kratt_1, Kratt_2} for the Hankel transform of a sequence with generating function of the form
$$\cfrac{1}{1-\alpha_1 x -\cfrac{ \beta_1 x^2}{1-\alpha_2 x - \cfrac{\beta_2 x^2}{1-\cdots}}} $$ can be derived from the above result, due to the fact that $ p_i =1$ in this case, and the fact that although in this note Equation (\ref{H}) has been expressed in the case of monomials $ b_i x^{p_i}$, the result continues to be true for polynomials $Q_{p_i}(x)=b_i x^{p_i}+\cdots$ of degree $p_i$.

\bigskip \hrule \bigskip

\noindent 2010 {\it Mathematics Subject Classification}: Primary 30B70; Secondary 11B83, 11C20, 15A15, 30B10
\\
\noindent \emph{Keywords:
C-fraction, continued fraction, power series, Hankel determinant, Hankel
transform, sequences}

 \bigskip \hrule \bigskip \noindent Concerns sequences

\seqnum{A000045},
\seqnum{A000108}

\end{document}